\documentstyle[12pt,leqno]{article}
\newcommand{\EE}{{\rm I\kern-2pt E}}
\newcommand{\RR}{{\rm I\kern-2pt R}}
\newcommand{\DD}{{\rm I\kern-2pt D}}
\newcommand{\PP}{{\rm I\kern-2pt P}}
\newcommand{\NN}{{\rm I\kern-2pt N}}
\newcommand{\dd}{{\rm \kern 3pt I\kern-9pt d}}

\newcommand{\Cc}{{\cal C}}
\newcommand{\Fc}{{\cal F}}
\newcommand{\Bc}{{\cal B}}

\newcommand{\Jc}{{\cal J}}

\title{CALCUL D'ERREUR COMPLET LIPSCHITZIEN ET FORMES DE DIRICHLET}
\author{Nicolas Bouleau\\ Ecole des Ponts, ParisTech}
\date{---}
\begin{document}
\maketitle

\vspace{.5cm}

\noindent{\bf Abstract}. We study the error calculus from a mathematical point
 of view, in particular for the infinite dimensional models
met in stochastic analysis, and also from the point of view of the links with experiments. Gauss was the first to propose an error calculus.
Thanks a coherence property this calculus is the most convenient in several questions. It can be reinforced by an extension
 principle based on Dirichlet forms which gives more strength to the coherence 
property. One gets a Lipschitzian complete error calculus which behaves well by images and by products and allows a quick and easy
construction of the basic mathematical tools of Malliavin's calculus. This allows also to catch the delicate question of error permanency 
that Poincar\'{e} emphasized. This error calculus is connected with statistics by mean of the notion of Fisher information. The article
ends with a comprehensive bibliography.

\begin{center}
---
\end{center}

\noindent{\bf R\'{e}sum\'{e}}. Nous \'{e}tudions le calcul d'erreur d'un point de vue math\'{e}matique en particulier pour les mod\`{e}les
de dimension infinie rencontr\'{e}s en calcul stochastique et \'{e}galement du point de vue du lien avec
 l'exp\'{e}rimentation. Gauss fut le premier \`{a} proposer un calcul
d'erreur au d\'{e}but du 19\`{e}me si\`{e}cle. Ce calcul poss\`{e}de une propri\'{e}t\'{e} de coh\'{e}rence qui le rend sup\'{e}rieur
dans bien des questions \`{a} d'autres formulations. Il peut \^{e}tre renforc\'{e} par un principe d'extension fond\'{e} sur la
 th\'{e}orie
des formes de Dirichlet. On obtient un calcul complet lipschitzien qui se comporte bien par image et par 
produit et permet une construction facile des notions de base du calcul de Malliavin. Cela permet d'aborder \'{e}galement la d\'{e}licate
question de la permanence des erreurs soulev\'{e}e par Poincar\'{e}. Ce calcul d'erreur se relie aux statistiques par l'interm\'{e}diaire
de l'information de Fisher. 

Le plan est le suivant : apr\`{e}s un aper\c{c}u des id\'{e}es de Gauss sur la loi des erreurs et de Poincar\'{e} sur la question de la 
permanence des erreurs et l'expos\'{e} du calcul de Gauss et de sa coh\'{e}rence, nous pr\'{e}sentons l'outil d'extension qui permet de 
construire un calcul lipschitzien et son axiomatisation. Nous abordons alors les liens avec l'exp\'{e}rimentation et les statistiques
puis les exemples d'\'{e}preuves r\'{e}p\'{e}t\'{e}es en dimension finie et infinie o\`{u} l'on rencontre le ph\'{e}nom\`{e}ne de permanence
des erreurs. Enfin nous \'{e}non\c{c}ons des conjectures induites par l'usage du calcul d'erreur lipschitzien pour d\'{e}montrer l'existence 
de densit\'{e}s sur l'espace de Wiener. L'article se termine par une bibliographie th\'{e}matique.\\

Au contraire des grandeurs discr\`{e}tes, les grandeurs continues sont le plus souvent entach\'{e}es d'erreur. Devant ce 
probl\`{e}me pratique plusieurs attitudes se rencontrent. Soit on traite les erreurs avec un langage vague, tra\c{c}ant
des barres d'erreur sans se soucier de corr\'{e}lation, ou prenant des maxima sans sp\'{e}cifier sur quel domaine, en justifiant
un tel laxisme du fait que les erreurs sont mal connues, soit on tente la gageure de propos rigoureux 
en d\'{e}gageant les hypoth\`{e}ses n\'{e}cessaires et les proc\'{e}dures d'exp\'{e}rimentation.

Cette seconde voie que nous allons suivre a \'{e}t\'{e} initi\'{e}e par Legendre, Laplace et Gauss au dŽ\'{e}but du 19\`{e}me
si\`{e}cle, dans une s\'{e}rie de travaux qu'on d\'{e}signe par {\it Th\'{e}orie classique des erreurs}. Le plus c\'{e}l\`{e}bre
d'entre eux est la d\'{e}monstration par Gauss de la ``loi des erreurs" par laquelle il montre, avec des hypoth\`{e}ses dont
certaines, implicites, seront relev\'{e}es par d'autres auteurs, que si l'on consid\`{e}re, dans une situation exp\'{e}rimentale,
que la moyenne arithm\'{e}tique des mesures faites est la meilleure valeur \`{a} prendre en compte, on doit admettre que 
les erreurs suivent une loi normale. Son raisonnement est probabiliste : la grandeur \`{a} mesurer est une variable al\'{e}atoire $X$ et
les mesures $X_1,\ldots, X_n$ sont suppos\'{e}es {\it conditionnellement ind\'{e}pendantes} sachant $X$.

A la fin du si\`{e}cle, dans son cours de {\it Calcul des probabilit\'{e}s} Henri Poincar\'{e} revient sur cette question en montrant que si 
on affaiblit certains pr\'{e}suppos\'{e}s de Gauss, d'autres lois que la lois normale peuvent \^{e}tre atteintes. Il discute longuement
un point nouveau et d\'{e}licat : le ph\'{e}nom\`{e}ne de {\it permanence des erreurs}. ``Avec un m\`{e}tre divis\'{e} en millim\`{e}tres,
on ne pourra jamais, \'{e}crit-il, si souvent qu'on r\'{e}p\`{e}te les mesures, d\'{e}terminer une longueur \`{a} un millioni\`{e}me de 
millim\`{e}tre pr\`{e}s". Ce ph\'{e}nom\`{e}ne est bien connu des physiciens, dans toute l'histoire de la physique on n'a jamais \'{e}t\'{e} 
capable de faire des mesures pr\'{e}cises avec des instruments grossiers cf [1]. Il ne d\'{e}veloppe pas de formalisme 
math\'{e}matique pour cela, il insiste en revanche sur l'avantage de supposer les erreurs petites car alors l'argument de Gauss devient
compatible avec les changements de variables non-lin\'{e}aires qui  peuvent s'\'{e}crire par le calcul diff\'{e}rentiel. C'est la question
du {\it calcul d'erreur}.\\

{\bf Le calcul d'erreur de Gauss}\\

Douze ans apr\`{e}s sa d\'{e}monstration conduisant \`{a} la loi normale, Gauss s'int\'{e}resse \`{a} la propagation des erreurs
({\it Theoria combinationis} 1821). Etant donn\'{e}e une grandeur $U=F(V_1,V_2,\ldots)$ fonction d'autres grandeurs $V_1,V_2,\ldots$, 
il pose le probl\`{e}me de calculer l'erreur quadratique de $U$ connaissant les erreurs quadratiques $\sigma_1^2,\sigma_2^2,\ldots$
de  $V_1,V_2,\ldots$ ces erreurs \'{e}tant suppos\'{e}es petites et ind\'{e}pendantes.

Sa r\'{e}ponse est la suivante
\begin{equation}
\sigma_U^2=(\frac{\partial U}{\partial V_1})^2\sigma_1^2 + 
(\frac{\partial U}{\partial V_2})^2\sigma_2^2 +
\ldots
\end{equation}
et il donne \'{e}galement la covariance de l'erreur de $U$ et d'une autre fonction des $V_1,V_2,\ldots.$

La formule (1) poss\`{e}de une propri\'{e}t\'{e} qui lui conf\`{e}re une grande sup\'{e}riorit\'{e} vis \`{a} vis d'autres formules souvent
propos\'{e}es dans les manuels. C'est la propri\'{e}t\'{e} de {\it coh\'{e}rence}. Avec une formule telle que
\begin{equation}
\sigma_U=|\frac{\partial U}{\partial V_1}|\sigma_1 + 
|\frac{\partial U}{\partial V_2}|\sigma_2 +
\ldots
\end{equation}
les erreurs peuvent d\'{e}pendre de la fa\c{c}on d'\'{e}crire la fonction $F$. En dimension 2 d\'{e}j\`{a} si on applique (2) \`{a} une application
lin\'{e}aire injective puis \`{a} son inverse, on obtient que l'identit\'{e} augmente les erreurs ce qui est difficilement
 acceptable.

Ceci ne se produit pas avec le calcul de Gauss. Pour le voir introduisons l'op\'{e}rateur diff\'{e}rentiel 
$$L=\frac{1}{2}\sigma_1^2\frac{\partial^2}{\partial V_1^2}+\frac{1}{2}\sigma_2^2\frac{\partial^2}{\partial V_2^2}+\ldots$$
et remarquons que (1) s'\'{e}crit
$$\sigma_U^2=LF^2-2FLF.$$
La coh\'{e}rence vient alors de la coh\'{e}rence du transport d'un op\'{e}rateur diff\'{e}rentiel par une fonction :  si $L$ est un tel 
op\'{e}rateur, si $u$ et $v$ d\'{e}signent des applications r\'{e}guli\`{e}res injectives et si on note $\theta_u L$ l'op\'{e}rateur 
$\varphi\rightarrow L(\varphi\circ u)\circ u^{-1}$, on a $\theta_{v\circ u}L=\theta_v(\theta_uL)$.

Les erreurs sur $V_1,V_2,\ldots$ peuvent ne pas \^{e}tre suppos\'{e}es ind\'{e}pendantes et peuvent d\'{e}pendre des valeurs de 
$V_1,V_2,\ldots$ : on se donne un champ de matrices sym\'{e}triques positives $(\sigma_{ij}(v_1,v_2,\ldots))$ sur $\RR^d$ repr\'{e}sentant
les variances et covariances conditionnelles de erreurs sachant les valeurs $v_1,v_2,\ldots$ de $V_1,V_2,\ldots$ et le calcul s'\'{e}crit

\begin{equation}
\sigma_F^2=\sum_{ij}\frac{\partial F}{\partial V_i}(v_1,v_2,\ldots)\frac{\partial F}{\partial V_j}(v_1,v_2,\ldots)\sigma_{ij}(v_1,v_2,\ldots)
\end{equation}
Le calcul d'erreur de Gauss traite des variances et covariances d'erreurs sans se pr\'{e}occuper des erreurs moyennes
 c'est-\`{a}-dire des biais. Il est important de souligner que {c'est la raison pour laquelle il ne fait intervenir que des d\'{e}riv\'{e}es 
premi\`{e}res}. En effet si on part d'une situation o\`{u} les erreurs sont centr\'{e}es, apr\`{e}s une application non lin\'{e}aire les erreurs
ne sont plus centr\'{e}es et le biais de l'erreur est du m\^{e}me ordre de grandeur que la variance. Par d'autres applications r\'{e}guli\`{e}res
non lin\'{e}aires cette situation va se perp\'{e}tuer. Ceci permet de voir que les variances peuvent se calculer par un calcul diff\'{e}rentiel
du premier ordre ne faisant intervenir que les variances, alors que les erreurs moyennes rel\`{e}vent d'un calcul du second ordre qui 
fait intervenir les moyennes et les variances.

La coh\'{e}rence du calcul de Gauss permet de le {\it g\'{e}om\'{e}triser}. Si une grandeur varie sur une vari\'{e}t\'{e} diff\'{e}rentiable
l'erreur associ\'{e}e peut \^{e}tre attach\'{e}e au point de la vari\'{e}t\'{e} comme objet g\'{e}om\'{e}trique. L'erreur
 est donn\'{e}e 
par une forme quadratique qui est une {\it m\'{e}trique riemannienne} sur la vari\'{e}t\'{e}. On peut prendre des {\it images}
par des applications injectives et de classe ${\cal C}^1$ en un calcul coh\'{e}rent ind\'{e}pendant des \'{e}critures des fonctions.
Ceci se relie \`{a} la th\'{e}orie des processus de diffusion sur les vari\'{e}t\'{e}s pour lesquelles nous
 renvoyons aux r\'{e}f\'{e}rences [2].\\

{\bf Calcul d'erreur avec outil d'extension}\\

Le calcul de Gauss est limit\'{e} par le fait qu'il ne dispose d'aucun moyen d'extension. A partir de l'erreur sur $(V_1,V_2,V_3)$
il permet de calculer l'erreur sur une fonction diff\'{e}rentiable de $(V_1,V_2,V_3)$ et c'est tout.

Or on aimerait \'{e}tendre ce calcul aux fonctions lipschitziennes
 car il est clair {\it a priori} qu'une application lipschitzienne de constante $\leq 1$, est contractante donc diminue les erreurs. Mais
surtout dans une situation fr\'{e}quente en calcul des probabilit\'{e}s o\`{u} on a une suite de quantit\'{e}s $X_1,X_2,\ldots,X_n,\ldots$
et o\`{u} on conna\^{\i}t les erreurs sur les fonctions r\'{e}guli\`{e}res d'un nombre fini de $X_n$, on aimerait pouvoir en d\'{e}duire l'erreur
sur des fonctions d'une infinit\'{e} des $X_n$ ou au moins sur certaines d'entre elles.

Il est en fait possible de doter le calcul d'erreur d'un puissant outil d'ex\-tension.

Pour cela on revient \`{a} l'id\'{e}e initiale de Gauss de consid\'{e}rer que les grandeurs \'{e}rron\'{e}es sont al\'{e}atoires, disons
d\'{e}finies sur $(\Omega,{\cal A}, \PP)$. L'erreur quadratique sur une variable al\'{e}atoire $X$ est elle-m\^{e}me al\'{e}atoire, nous
la notons $\Gamma[X]$. Elle est infinit\'{e}simale mais cela n'appara\^{\i}t pas dans les notations, comme si nous avions
une unit\'{e} de mesure infinit\'{e}simale pour les erreurs fix\'{e}e dans tout le probl\`{e}me. L'outil est le suivant : nous supposons
que si $X_n\rightarrow X$ dans $L^2(\Omega,{\cal A}, \PP)$ et si l'erreur $\Gamma[X_m-X_n]$ sur $X_m-X_n$ peut \^{e}tre rendue aussi petite
qu'on veut dans $L^1$ pour $m,n$ grands, alors l'erreur $\Gamma[X_n-X]$ tend vers z\'{e}ro dans $L^1$.

C'est un {\it principe de coh\'{e}rence renforc\'{e}} puisqu'il signifie que l'erreur quadratique sur $X$ est attach\'{e}e \`{a} $X$ en tant qu'application
math\'{e}matique et que si le couple $(X_n,{\mbox{ erreur quadratique sur }}X_n)$ converge en un sens convenable il converge nŽcessairement vers 
$(X,{\mbox{ erreur quadratique sur }}X)$.

Ceci s'axiomatise de la fa\c{c}on suivante.

Nous appellerons {\it structure d'erreur} un espace de probabilit\'{e} muni
 d'une forme de Dirichlet locale poss\'{e}dant un op\'{e}rateur carr\'{e} du champ.
Plus pr\'{e}ci\-s\'{e}ment c'est un terme
$$(\Omega,\,{\cal{A}},\,\PP,\,\DD,\,\Gamma)$$ o\`{u} $(\Omega,{\cal{A}},\PP)$
 est un espace de probabilit\'{e}, v\'{e}rifiant les quatre propri\'{e}t\'{e}s suivantes :

1. $\DD$ est un sous-espace vectoriel dense de $L^2(\Omega,{\cal{A}},\PP)$

2. $\Gamma$ est une application bilin\'{e}aire sym\'{e}trique positive de $\DD\times\DD$
 dans $L^1(\PP)$ v\'{e}rifiant le calcul fonctionnel de classe ${\cal{C}}^1\cap{\mbox{Lip}}$,
ce qui signifie que si $u\in\DD^m$ et $v\in\DD^n$ pour $F$ et $G$ de classe ${\cal{C}}^1$ et lipschitziennes de 
$\RR^m$ [resp. $\RR^n$] dans $\RR$, on a $F\circ u\in\DD$ et $G\circ v\in\DD$ et
$$\Gamma[F\circ u,G\circ v]=\sum_{i,j} F_i^{\prime}(u) G_j^{\prime}(v) \Gamma[u_i,v_j]\quad\PP{\mbox{-p.s.}}$$

3. La forme bilin\'{e}aire ${\cal{E}}[f,g]=\EE[\Gamma[f,g]]$ est ferm\'{e}e, ce qui signifie que 
$\DD$ est complet pour la norme $\|\,.\,\|_{\DD}=(\|\,.\,\|_{L^2(\PP)}^2 +{\cal{E}}[\,.\,,\,.\,])^{\frac{1}{2}}$.

4. $1\in\DD$ et $\Gamma[1,1]=0$.\\

On notera ${\cal{E}}[f]$ pour ${\cal{E}}[f,f]$ et $\Gamma[f]$ pour $\Gamma[f,f]$.\\

\noindent{\it Commentaire.}  Avec cette d\'{e}finition la forme ${\cal{E}}$ est une {\it forme
 de Dirichlet} notion  introduite par Beurling et Deny comme outil de th\'{e}orie du potentiel et qui re\c{c}ut
une interpr\'{e}tation probabiliste en termes de processus de Markov sym\'{e}trique par les 
travaux de Silverstein et Fukushima cf [3]. L'op\'{e}rateur $\Gamma$ est le carr\'{e} du champ
 associ\'{e} \`{a} ${\cal{E}}$, \'{e}tudi\'{e} par de nombreux auteurs dans des contextes plus g\'{e}n\'{e}raux que celui-ci cf [4]. Les formes de Dirichlet locales et leur
carr\'{e} du champ admettent un calcul fonctionnel de classe ${\cal{C}}^1\cap{\mbox{Lip}}$ dont l'interpr\'{e}tation
probabiliste sort du cadre des semi-martingales cf [5].\\

\noindent{\it Premiers exemples}. a)  Un exemple simple de structure d'erreur est le terme $$(\RR,\,{\cal{B}}(\RR),\,\mu,\,H^1(\mu),\,\gamma)$$
 o\`{u} $\mu=N(0,1)$ et $H^1(\mu)=\{f\in L^2(\mu)\,:\,f^{\prime}{\mbox{ (au sens de }}{\cal{D}}^{\prime})\in L^2(\mu)\}$ 
avec $\gamma[f]=f^{\prime 2}$ pour $f\in H^1(\mu)$. Cette structure est associ\'{e}e au processus d'Ornstein-Uhlenbeck
 \`{a} valeurs r\'{e}elles.

b) Soit $D$ un ouvert connexe de $\RR^d$ de volume unit\'{e}, $\lambda_d$ la mesure de Lebesgue, on prend
$(\Omega,{\cal A},\PP)=(D,{\cal B}(D),\lambda_d)$. On d\'{e}finit
$$\Gamma[u,v]=\sum_{ij}\frac{\partial u}{\partial x_i}
\frac{\partial v}{\partial x_j}a_{ij}\quad\quad{\mbox{pour }}u,v\in{\cal C}_K^{\infty}(D)$$
o\`{u} les $a_{ij}$ sont des  applications de $D$ dans $\RR$ telles que 
$$a_{ij}\in L_{loc}^2(D),\; a_{ij}=a_{ji},\;\frac{\partial a_{ij}}{\partial x_k}\in L_{loc}^2(D),
\;\sum_{ij}a_{ij}(x)\xi_i\xi_j\geq 0\;\forall \xi\in\RR^d
\; \forall x\in D.$$
On peut alors montrer que la forme ${\cal E}[u,v]=\EE\Gamma[u,v]$ avec $u,v\in{\cal C}_K^{\infty}(D)$ est fermable (cf [6]) autrement dit il
existe une extension de $\Gamma$ \`{a} un sous-espace $\DD$ de $L^2$, $\DD\supset{\cal C}_K^\infty(D)$ telle que 
$(\Omega,\,{\cal{A}},\,\PP,\,\DD,\,\Gamma)$ soit une structure d'erreur.\\

C'est une cons\'{e}quence des hypoth\`{e}ses des structures d'erreur que si $f\in\DD$ et si $F$ est
 lipschitzienne de $\RR$ dans $\RR$ alors $F\circ f\in\DD$ et $\Gamma[F\circ f]\leq k\Gamma[f]$. Plus g\'{e}n\'{e}ralement
 si $F$ est une contraction de $\RR^d$ dans $\RR$ au sens suivant
$$|F(x)-F(y)|\leq \sum_{i=1}^d|x_i-y_i|$$
alors si $f_1,f_2,\ldots,f_d\in\DD$ on a $F(f_1,f_2,\ldots,f_d)\in\DD$ et
$$\Gamma[F(f_1,f_2,\ldots,f_d)]^{\frac{1}{2}}\leq\sum_{i=1}^d\Gamma[f_i]^{\frac{1}{2}}.$$

Deux propri\'{e}t\'{e}s facilitent le maniement des structures d'erreur et permettent d'accompagner les constructions de 
mod\`{e}les probabilistes : 

1) L'op\'{e}ration de {\it prendre l'image} d'une structure d'erreur par une application se fait tr\`{e}s naturellement
et fournit encore une structure d'erreur, d\`{e}s que l'application satisfait certaines
 conditions assez larges cf [7]. En particulier si $$(\Omega,\,{\cal{A}},\,\PP,\,\DD,\,\Gamma)$$ 
est une structure d'erreur et si $X$ est une variable al\'{e}atoire \`{a} valeurs $\RR^d$ dont les composantes sont dans $\DD$
$$(\RR^d,\,{\cal{B}}(\RR^d),\,\PP_X,\,\DD_X,\,\Gamma_X)$$ est une structure d'erreur o\`{u} $\PP_X$ est la
loi de $X$, 
$$\begin{array}{rcl}
\DD_X&=&\{f\in L^2(\PP_X)\,:\,f\circ X\in\DD\}\\
\Gamma_X[f](x)&=&\EE[\Gamma[f\circ X]|X=x],\quad\quad f\in\DD_X.
\end{array}
$$ 

Soulignons que dans ce calcul d'erreur complet lipschitzien les images ne sont pas limit\'{e}es \`{a} des applications injectives. Par 
exemple la structure de l'exemple a) $(\RR,\,{\cal{B}}(\RR),\,\mu,\,H^1(\mu),\,\gamma)$ a une image par 
l'application $$x\rightarrow|\sin{\sqrt{1+|x|}}|$$ qui est une structure d'erreur sur $[0,1]$.\\

2) {\it Le produit} de deux ou d'une infinit\'{e} d\'{e}nombrable de structures d'erreur est toujours d\'{e}fini et donne une structure
 d'erreur. On obtient ainsi facilement des structures d'erreurs sur des espaces de dimension
 infinie, cf [7], par exemple sur l'espace de Wiener ou sur l'espace
 de Poisson et sur les mod\`{e}les qui s'en d\'{e}duisent, c'est une fa\c{c}on d'aborder
 le calcul de Malliavin, cf [8]. 

Indiquons \`{a} titre d'exemple la construction de la structure d'Ornstein-Uhlenbeck sur l'espace de Wiener.

Reprenons la structure d'erreur unidimensionnelle de l'exemple a)  
$$ (\RR,{\cal B}(\RR),\mu, H^1(\mu),\gamma)$$ et consid\'{e}rons la structure produit infini qui s'en d\'{e}duit : 
$$(\Omega,{\cal A}, \PP, \DD,\Gamma)=(\RR,{\cal B}(\RR), \mu, H^1(\mu),\gamma)^{\NN}
=(\RR^{\NN},{\cal B}(\RR^{\NN}), \mu^{\NN}, \DD, \Gamma).$$
 Les applications coordonn\'{e}es $X_n$, par construction du produit,  sont gaussien\-nes r\'{e}duites ind\'{e}pendantes,
appartiennent \`{a} $\DD$ et v\'{e}rifient
$$\begin{array}{rcl}
\Gamma[X_n]&=&1\\
\Gamma[X_m,X_n]&=&0\quad m\neq n
\end{array}
$$
Soit $\chi_n$une base orthonormale de $L^2(\RR_+,dt)$. On pose
$$B_t=\sum_{n=0}^\infty\int_0^t\chi_n(s)\,ds.X_n.$$
$(B_t)_{t\geq 0}$ est un mouvement brownien et si  $f\in L^2(\RR_+)$ s'\'{e}crit $f=\sum_n a_n\chi_n$, la variable al\'{e}atoire
$\sum_na_nX_n$ est not\'{e}e $\int_0^\infty f(s)\,dB_s$ par extension du cas o\`{u} $f$ est \'{e}tag\'{e}e.

Nous avons alors $\int f(s)\,dB_s\in\DD$ et
$$\Gamma[\int f(s)\,dB_s]=\Gamma[\sum_na_nX_n]=\sum_na_n^2\Gamma[X_n]=\sum_n a_n^2=
\|f\|_{L^2(\RR_+,dt)}^2$$
Puis par les r\`{e}gles du calcul d'erreur si $F\in{\cal C}^1\cap{\mbox{Lip}}(\RR^d)$
$$\Gamma(F(\int\! f_1(s)dB_s,\ldots,\int\! f_d(s)dB_s)]=
\sum_i F_i^{\prime 2}(\int\! f_1(s)dB_s,\ldots,\int\! f_d(s)dB_s)\|f_i\|_{L^2(dt)}^2$$
et par l'outil d'extension, le calcul d'erreur s'\'{e}tend \`{a} d'autres fonctionnelles browniennes dont les solutions
 d'\'{e}quations diff\'{e}rentielles stochastiques \`{a} coefficients lipschitziens cf [7] [8].\\

{\bf Calculs d'erreur et statistiques}\\

Pour passer du calcul d'erreur de Gauss au calcul complet lipschitzien il est n\'{e}cessaire de disposer d'une probabilit\'{e}.
Si des grandeurs sont variables mais d\'{e}terministes comme parfois en m\'{e}canique, elles doivent \^{e}tre replac\'{e}es dans un cadre
probabiliste. C'est le terme $(\Omega,{\cal A},\PP)$ d'une structure d'erreur $(\Omega,{\cal A},\PP,\DD,\Gamma)$.

Une premi\`{e}re approche consiste \`{a} suivre les id\'{e}es de E. Hopf dans les ann\'{e}es 1930 qui dans l'esprit des travaux de Poincar\'{e}
montra, par des formes g\'{e}n\'{e}rales de th\'{e}or\`{e}mes limites en loi, que de nombreux syst\`{e}mes dynamiques poss\`{e}dent des 
lois de probabilit\'{e} naturelles qu'on peut prendre comme loi {\it a priori}, cf [9].

Une seconde voie consiste \`{a} se donner un op\'{e}rateur elliptique du second ordre $L$ v\'{e}rifiant
\begin{equation}
\Gamma[F]=LF^2-2FLF
\end{equation}
(dont seuls les termes du second ordre sont d\'{e}termin\'{e}s par cette relation) qui fournit le cadre d'un calcul d'erreur pour les 
variances {\it et les biais}. Puis de construire la probabilit\'{e} invariante vis \`{a} vis de laquelle la diffusion de g\'{e}n\'{e}rateur
$L$ est sym\'{e}trique. Ceci peut \^{e}tre fait de fa\c{c}on g\'{e}om\'{e}trique en dimension
 finie ou infinie cf [10].

Nous suivons une troisi\`{e}me voie qui se relie plus directement aux applications. Nous consid\'{e}rons que les conditions exp\'{e}rimentales
sont suffisamment sp\'{e}cifi\'{e}es pour que la probabilit\'{e} $\PP$ s'obtienne comme habituellement par les statistiques et
nous allons montrer que les statistiques fournissent en fait \'{e}galement l'op\'{e}rateur $\Gamma$ donc finalement la structure 
d'erreur, au moins sur un domaine minimal pour $\Gamma$.

Consid\'{e}rons une grandeur erron\'{e}e $d$-dimensionnelle $X$. L'espace image par $X$ est 
$$(\RR^d,\Bc(\RR^d),\PP_X(dx))$$
L'op\'{e}rateur $\Gamma$ que nous cherchons \`{a} d\'{e}finir se pr\'{e}sente sous la forme
$$\Gamma_X[F](x)=\sum_{i,j=1}^d F_i^{\prime}(x)F_j^{\prime}(x)a_{ij}(x)$$
o\`{u} la matrice $A(x)=(a_{ij}(x))$ est sym\'{e}trique positive, c'est elle qu'il faut conna\^{\i}tre
et qui repr\'{e}sente la pr\'{e}cision avec laquelle $X$ est connu au point $x$.

Remarquons que si $G\,:\,\RR^d\mapsto\RR^m$ est $\Cc^1\cap{\mbox{Lip}}$, d'apr\`{e}s le calcul fonctionnel,
 la variable al\'{e}atoire $G(x)$ est alors connue avec la pr\'{e}cision\\

\noindent$(5)\hfill \Gamma_X[G, G^t](x)=\nabla_xG.A(x).(\nabla_xG)^t\hspace{4cm}$\\

\noindent o\`{u} $\nabla_xG$ est la matrice jacobienne de $G$ en $x$. 

Mais pour conna\^{\i}tre $X$, sous la loi conditionnelle $X=x$ not\'{e}e $\EE_x$, nous proc\'{e}dons \`{a}
 des mesures qui sont {\it des estimateurs} du param\`{e}tre $x$. Soit $T$ un tel
estimateur \`{a} valeur $\RR^m$ de matrice de covariance $\EE_x[(T-\EE_x[T]). (T-\EE_x[T])^t]$. Sous les
 hypoth\`{e}ses statistiques dites du {\it mod\`{e}le r\'{e}gulier} l'in\'{e}galit\'{e} de 
Fr\'{e}chet-Darmois-Cramer-Rao s'\'{e}crit\\

\noindent$(6)\hfill\EE_x[(T-\EE_x[T]). (T-\EE_x[T])^t]\geq \nabla_x\EE_xT.J(x)^{-1}.(\nabla_x\EE_xT)^t\hspace{1cm}$\\

\noindent au sens de l'ordre du c\^{o}ne des matrices sym\'{e}triques positives, o\`{u} $J(x)$ est la matrice
 d'{\it information de Fisher}, cf [11]. La  meilleure pr\'{e}cision qu'on peut avoir sur
$X$ est donc $J(x)^{-1}$ et la comparaison de 
(5) et (6) conduit \`{a} poser
$$A(x)=J(x)^{-1}.$$ On se convainc facilement que cette d\'{e}finition est compatible avec les changements de
 variables r\'{e}guliers
 : si on estime $\psi(x)$ au lieu de $x$, on obtient comme structure d'erreur l'image par $\psi$ de la structure d'erreur de $X$.

Cette connexion entre l'information de Fisher et   l'approche des erreurs fond\'{e}e sur les formes de Dirichlet
pose une s\'{e}rie de questions qui sont encore au stade de la recherche. J'en \'{e}noncerai
 trois :

a) Sous quelles hypoth\`{e}ses peut-on obtenir directement $J(x)^{-1}$ \'{e}ventuel\-lement singuli\`{e}re sans passer
 par la matrice d'information de Fisher ?

b) Les m\'{e}thodes de statistique asymptotique cf [12] donnent-elles des outils pour \'{e}tudier
 la fermabilit\'{e} des pr\'{e}-formes de Dirichlet sur $\RR^d$ ?

c) L'emploi en statistique de mod\`{e}les qui sont des structures d'erreur (ferm\'{e}es par hypoth\`{e}se)
permet-il de pr\'{e}ciser certains th\'{e}or\`{e}mes asymptotiques ?\\

{\bf Syst\`{e}mes projectifs et \'{e}preuves r\'{e}p\'{e}t\'{e}es}\\

L'introduction d'op\'{e}rateurs d'erreur en plus du langage probabiliste permet de traiter avec beaucoup de finesse la question des \'{e}preuves
r\'{e}p\'{e}t\'{e}es et de r\'{e}pondre par des mod\'{e}lisations explicites au ph\'{e}nom\`{e}ne de permanence des erreurs point\'{e} par 
Poincar\'{e}.

Comme nous allons le voir certains syst\`{e}mes projectifs pour lesquels la limite projective des espaces
 de proba\-bilit\'{e} existe, n'admettent pas de structure d'erreur limite, mais d\'{e}finissent
seulement une pr\'{e}-structure d'erreur au sens suivant :

Un terme $(\Omega,\,{\cal{A}},\,\PP,\,\DD^0,\,\Gamma)$ est une pr\'{e}-structure d'erreur si
$(\Omega,{\cal{A}},\PP)$ est une espace de probabilit\'{e} et si $\DD^0,\,\Gamma$ v\'{e}rifient les propri\'{e}t\'{e}s 
(1.), (2.) et (4.) des structures d'erreurs mais pas n\'{e}cessairement la propri\'{e}t\'{e} (3.).

Il y a donc des pr\'{e}-structures d'erreur fermables et des pr\'{e}-structures d'erreur non-fermables. Images
 et produits se d\'{e}finissent facilement pour des pr\'{e}-structures d'erreur.

Fixons quelques notations pour les syst\`{e}mes projectifs sous les hypoth\`{e}ses
 de r\'{e}gularit\'{e} courantes cf [7] : 

Etant donn\'{e}s des espaces mesurables $(E_i,{\cal{F}}_i)\;i\in \NN^\ast$, pour
 $\alpha\in{\cal{J}}$
ensemble des parties finies de $\NN^\ast$ un syst\`{e}me projectif de structures
 d'erreur (ou de pr\'{e}-structures d'erreur) est une famille 
$$(E_{\alpha},\,{\cal{F}}_{\alpha},\,m_{\alpha},\,\DD^0_{\alpha},\,\Gamma_{\alpha})$$ de 
(pr\'{e}-)structures d'erreur o\`{u} $(E_{\alpha},{\cal{F}}_{\alpha})=
\prod_{i\in\alpha}(E_i,{\cal{F}}_i)$
qui sont compatibles au sens usuel. Posant 
$$\DD^0=\cup_{\alpha\in{\cal{J}}}\DD^0_{\alpha}$$ il d\'{e}finit
une pr\'{e}-structure d'erreur
$$(E,\,{\cal{F}},\,m,\,\DD^0,\,\Gamma)$$ dont les projections sont
 les $(E_{\alpha},{\cal{F}}_{\alpha},m_{\alpha},\DD^0_{\alpha},\Gamma_{\alpha})$.

Les syst\`{e}mes projectifs que nous consid\'{e}rons par la suite sont tels que les 
$(E_i,{\cal{F}}_i)$ sont identiques et que le syst\`{e}me projectif soit auto-isomorphe par translation des indices,
ceci afin de repr\'{e}senter des \'{e}preuves r\'{e}p\'{e}t\'{e}es.

Nous allons donner trois exemples. 
Dans les exemples A et B la situation dont on prend des \'{e}preuves r\'{e}p\'{e}t\'{e}es est finie dimensionnelle, c'est un 
mod\`{e}le probabiliste de dimension finie avec des grandeurs \'{e}rron\'{e}es. Les propri\'{e}t\'{e}s
asymptotiques des \'{e}preuves r\'{e}p\'{e}t\'{e}es sont diff\'{e}rentes dans les cas A et B. Dans l'exemple
 C le mod\`{e}le probabiliste dont on fait des \'{e}preuves r\'{e}p\'{e}t\'{e}es est un espace de processus al\'{e}atoire,
 infini-dimensionnel, ce cas est important car il donne l'id\'{e}e des applications les plus int\'{e}ressantes
 (physique statistique, filtrage et pr\'{e}diction, finance).\\

\noindent{\bf A}. Dans le premier exemple que nous prenons, les erreurs sont corr\'{e}l\'{e}es, (ainsi que le sugg\'{e}rait
 Poincar\'{e}) et le syst\`{e}me projectif est fermable :

$$(E_i,{\cal{F}}_i)=([0,1],{\cal{B}}([0,1]))\quad\forall i\in\NN^{\ast}$$.

\noindent La pr\'{e}-structure $(E_{\alpha},{\cal{F}}_{\alpha}
,m_{\alpha},\DD^0_{\alpha},\Gamma_{\alpha})$
est ainsi d\'{e}finie
$$(E_{\alpha},{\cal{F}}_{\alpha},m_{\alpha})=([0,1]^{|\alpha|},
{\cal{B}}([0,1])^{|\alpha|},\lambda_{|\alpha|})$$
o\`{u} $|\alpha|={\mbox{card}}(\alpha)$ et $\lambda_{|\alpha|}$ est
 la mesure de Lebesgue de dimension $|\alpha|$,

$$\DD_{\alpha}^0={\cal{C}}_K^0(]0,1[^{|\alpha|})\oplus\RR$$

\noindent et pour $u,v\in\DD_{\alpha}^0$

$$\Gamma_{\alpha}[u,v]=\sum_{i,j\in\alpha}\frac{\partial u}{\partial x_i}\frac{\partial v}{\partial x_j}a_{ij}$$

\noindent o\`{u} les $a_{ij}$ sont constants et tels que les matrices $(a_{ij})_{i,j\in\alpha}$
 soient sym\'{e}triques  positives.

Dans ce cas on peut montrer que les pr\'{e}-structures  $(E_{\alpha},{\cal{F}}_{\alpha},m_{\alpha},\DD^0_{\alpha},\Gamma_{\alpha})$
sont fermables et \'{e}galement la pr\'{e}-structure d\'{e}finie par le syst\`{e}me projectif 
$$(E,\,{\cal{F}},\,m,\,\DD^0,\,\Gamma)$$
sa fermeture d\'{e}finit une structure d'erreur
$$(E,\,{\cal{F}},\,m,\,\DD,\,\Gamma)$$
 dont les projections sont les fermetures des 
$(E_{\alpha},{\cal{F}}_{\alpha},m_{\alpha},\DD^0_{\alpha},\Gamma_{\alpha})$. La d\'{e}monstration est fond\'{e}e sur le fait que
$$Au=\sum_{ij}\frac{\partial}{\partial x_i}(a_{ij}\frac{\partial u}{\partial x_j})\quad u\in\DD^0$$
d\'{e}finit un op\'{e}rateur sym\'{e}trique ce qui, par un argument classique cf [5], donne le r\'{e}sultat.

Examinons-la structure d'erreur $(E,\,{\cal{F}},\,m,\,\DD,\,\Gamma)$ plus en d\'{e}tail.

Soient $(U_n)_{n\geq 1}$ les applications coordonn\'{e}es. Les $U_n$ ne sont pas dans $\DD$
 mais si $\psi\in H_0^1(]0,1[)$ les  
variables al\'{e}atoires $X_n=\psi(U_n)$ sont dans $\DD$, elles sont i.i.d. et
$$\begin{array}{rcl}
\Gamma[X_n]&=&\psi^{\prime 2}(U_n)a_{nn}\\
\Gamma[X_m,X_n]&=&\psi^{\prime}(U_m)\psi^{\prime}(U_n)a_{mn}.
\end{array}
$$
a) Si $a_{ij}=1\quad\forall i,j$ nous avons $\forall F\in{\cal{C}}^1\cap{\mbox{Lip}}(\RR^d)$
$$\Gamma[F(X_1,\ldots,X_d)]=\left(\sum_{i=1}^d F_i^{\prime}(X_1,\ldots,X_d)\psi^{\prime}(U_i)\right)^2$$
ainsi

$$\lim_{N\uparrow\infty}\frac{1}{N}\sum_{n=1}^N X_n=\int_0^1 \psi(x)\,dx$$

\noindent et d'autre part

$$\lim_{N\uparrow\infty}\Gamma[\frac{1}{N}\sum_{n=1}^N X_n]=\lim_{N\uparrow\infty}(\frac{1}{N}
\sum_{n=1}^N\psi^{\prime}(U_n))^2=(\int_0^1\psi^{\prime}(x)\,dx)^2=0$$

\noindent les limites \'{e}tant dans $L^1$ et p.s. par la loi des grands nombres.

\noindent b) Plus g\'{e}n\'{e}ralement dans le cas $a_{ij}=a(i-j)$ avec $\sum \xi_i\xi_j a(i-j)\geq 0$ utilisant le th\'{e}or\`{e}me de 
repr\'{e}sentation de Bochner, on obtient aussi

$$\lim_{N\uparrow\infty}\Gamma[\frac{1}{N}\sum_{n=1}^N X_n]=0 {\mbox{ dans }} L^1$$.

On voit que dans cet exemple les \'{e}preuves sont ind\'{e}pendantes avec des erreurs corr\'{e}l\'{e}es
 mais l'erreur sur la moyenne s'\'{e}vanouit. Ce ne sera plus le cas dans l'exemple suivant.\\

\noindent {\bf B}. Supposons que chaque tirage concerne une quantit\'{e} scalaire et soit
$$(\RR,{\cal{B}}(\RR),\rho,\dd,G)$$ une structure d'erreur sur $\RR$ telle que l'identit\'{e} soit
dans $\dd$ et v\'{e}rifiant pour $u\in\Cc^1\cap{\mbox{Lip}}\subset\dd$
$$G[u]=u^{\prime 2}.g$$ o\`{u} $g$ est une fonction positive de $L^1(\rho)$.

D\'{e}finissons la pr\'{e}-structure $(E_{\alpha},{\cal{F}}_{\alpha},m_{\alpha},\DD^0_{\alpha},\Gamma_{\alpha})$
par 

$(E_{\alpha},{\cal{F}}_{\alpha},m_{\alpha})=(\RR,{\cal{B}}(\RR),\rho)^{|\alpha|}$,

\noindent soient $(X_n)_{n\geq 1}$ les applications coordonn\'{e}es, prenons

$\DD_{\alpha}^0=\Cc^1\cap{\mbox{Lip}}(\RR^{|\alpha|})$

\noindent et pour $u\in\DD_{\alpha}^0$ posons

$\Gamma{_\alpha}[u]=(\sum_{i\in\alpha} u_i^{\prime} f(X_i))^2+\sum_{i\in\alpha}u_i^{\prime 2}.g(X_i)
\quad{\mbox{pour une fonction }} f\in L^2(\rho).$

On peut montrer, par exemple en supposant les fonctions $f,g$ telles que $\frac{f^2}{g}$ soit born\'{e}e,
que les pr\'{e}-structures $(E_{\alpha},{\cal{F}}_{\alpha},m_{\alpha},\DD^0_{\alpha},\Gamma_{\alpha})$
sont fermables et leurs fermetures $(E_{\alpha},{\cal{F}}_{\alpha},m_{\alpha},\DD_{\alpha},\Gamma_{\alpha})$
forment un syst\`{e}me projectif d\'{e}finissant une pr\'{e}-structure
$$(E,\Fc,m,\DD^0,\Gamma)\quad\quad\DD^0=\cup_{\alpha\in\Jc}\DD_{\alpha}.$$ 
Les coordonn\'{e}es $X_n$ sont dans $\DD^0$ et sont i.i.d. de loi $\rho$.

Nous avons 
$$\begin{array}{rcl}
\Gamma[X_n]&=&f^2(X_n)+g(X_n)\\
\Gamma[X_m,X_n]&=&f(X_m)f(X_n)\quad m\neq n.
\end{array}
$$
 Si $h\in\Cc^1\cap{\mbox{Lip}}(\RR)$ on a lorsque $N\uparrow\infty$
$$\begin{array}{rcll}
\frac{1}{N}\sum_{n=1}^N h(X_n)&\rightarrow& \int h\,d\rho&{\mbox{p.s. et dans }}L^2(m)\\
\Gamma[\frac{1}{N}\sum_{n=1}^N h(X_n)]&\rightarrow&(\int h^{\prime} .f.d\rho)^2& {\mbox{ p.s. et dans }} L^1(m).
\end{array}
$$
Or si $u_N=\frac{1}{N}\sum_{n=1}^N h(X_n)$ on peut montrer ais\'{e}ment que $\Gamma[u_M-u_N]\rightarrow 0$ dans $L^1$
 pour $M,N\uparrow\infty$. Il en r\'{e}sulte que la pr\'{e}-structure $(E,\Fc,m,\DD^0,\Gamma)$ {\it n'est jamais fermable} \`{a}
moins que $f=0\quad{\mbox{p.s. }}$ ce qui redonne le cas d'une structure produit.

 La limite $(\int h^{\prime} .f.d\rho)^2$ n'est pas nulle en g\'{e}n\'{e}ral, le mod\`{e}le d\'{e}crit une situation
 analogue \`{a}
 celle relev\'{e}e par Poincar\'{e} o\`{u} les erreurs ne sont pas \'{e}vanescentes.\\

\noindent{\bf C}. Passons \`{a} des exemples  faisant intervenir la dimension infinie.

Une ficelle de longueur $L$ est jet\'{e}e sur le plan et on mesure (par exemple au moyen de lignes
parall\`{e}les tr\`{e}s resserr\'{e}es) la longueur totale de la projection de la ficelle sur $Ox$.
Par \'{e}preuves r\'{e}p\'{e}t\'{e}es ceci permet de mesurer la longueur de la ficelle comme nous allons le voir dans un instant.

Ceci peut \^{e}tre mod\'{e}lis\'{e} de la fa\c{c}on suivante :

\noindent la ficelle est param\'{e}tr\'{e}e par 
$$\begin{array}{rcl}
X(t)&=&X_0+\int_0^t\cos(\varphi+B_s)\,ds\\
Y(t)&=&Y_0+\int_0^t\sin(\varphi+B_s)\,ds\quad\quad 0\leq t\leq L\leq 1
\end{array}
$$ o\`{u} $B$ est un mouvement brownien standard et $\varphi$ uniforme sur le cercle,
 ind\'{e}pendant de $B$. On mesure la quantit\'{e}
$$A(\varphi,\omega)=\int_0^L|\cos(\varphi+B_s)|\,ds$$ On obtient $\EE A$ par \'{e}preuves r\'{e}p\'{e}t\'{e}es
 et on en d\'{e}duit la longueur $L$ de la ficelle par la formule $$\EE A=\frac{2L}{\pi}$$ qui
vient imm\'{e}diatement de l'expression de $A$ par int\'{e}gration puisque $\varphi$ et $B$ sont ind\'{e}pendants.

Comme hypoth\`{e}se sur les erreurs nous supposons qu'il y a une erreur sur $\varphi$ et
 une erreur sur $B$ ind\'{e}pendantes mais que les erreurs sur les diverses \'{e}preuves sont corr\'{e}l\'{e}es.

Sur $\varphi$ on consid\`{e}re une erreur analogue au cas B.

Sur $B$ on consid\`{e}re pour la simplicit\'{e} l'erreur donn\'{e}e par la forme de Dirichlet
associ\'{e}e au semigroupe d'Ornstein-Uhlenbeck cf [8] dont nous avons indiqu\'{e} la 
construction plus haut.

Notons $(\Omega,\Fc,\PP)$ l'espace de Wiener. Cette structure d'erreur  est telle que pour une variable al\'{e}atoire
 de la forme $\int_0^1 f(s)\,dB_s$, $\,f\!\in \!L^2([0,1]),$
l'op\'{e}rateur $\Gamma_{OU}$ ($OU$ pour Ornstein-Uhlenbeck) est d\'{e}fini par 
$$\Gamma_{OU}[\int_0^1 f(s)\,dB_s]=\int_0^1 f^2(s)\,ds$$
puis la d\'{e}finition de $\Gamma_{OU}$ s'\'{e}tend par le calcul fonctionnel aux variables de
 la forme
$$F(\int_0^1 f_1(s)\,dB_s,\int_0^1 f_2(s)\,dB_s,\ldots, \int_0^1 f_k(s)\,dB_s)$$ et \`{a} d'autres
 par continuit\'{e} de $\Gamma_{OU}$
dont le domaine est not\'{e} $\DD_2^1$. On peut aussi construire un op\'{e}rateur gradient $D$ d\'{e}fini
sur $\DD_2^1$ \`{a} valeurs dans $L^2(\PP,H)$ o\`{u} $H$ est un espace de Hilbert auxiliaire
 et reli\'{e} \`{a} $\Gamma_{OU}$
par $$\forall G\in\DD_2^1\quad\Gamma_{OU}[G]=<DG,DG>_H.$$
Nous notons $\varphi_n,\omega_n$ les applications coordonn\'{e}es du syst\`{e}me projectif
 et si $X(\varphi,\omega)$ est
 une variable al\'{e}atoire d\'{e}pendant du lacet $(\varphi,\omega)$ nous notons
$$X_n=X(\varphi_n,\omega_n)$$
le $n$-i\`{e}me tirage de cette variable al\'{e}atoire.

La corr\'{e}lation des erreurs peut \^{e}tre exprim\'{e}e au moyen d'un op\'{e}rateur de Hilbert-Schmidt $K$ sur $H$
d\'{e}pendant de $\omega$ tel que $\|K\|_{HS}\in L^\infty.$

On consid\`{e}re $\prod_{n=1}^N(\Omega,\Fc,\PP)_n$ (copies de l'espace de Wiener) et pour une variable al\'{e}atoire $F(\omega_1,\ldots,\omega_N)$
nous d\'{e}finissons l'op\'{e}rateur $\Gamma_N$ par 
$$\Gamma_N[F]=\|\sum_{i=1}^N K_iD_iF\|_H^2+\sum_{i=1}^N\|D_iF\|_H^2$$
o\`{u} $K_i$ et $D_i$ sont les op\'{e}rateurs $K$ et $D$ op\'{e}rant sur $\omega_i$. Soit\\

\noindent$\DD_N=\{F\in L^2(\PP^N)\;:\;\forall i\;\omega_i\mapsto F(\omega_1,\ldots,\omega_i,\ldots,\omega_N)\in\DD_2^1\hspace*{\fill}$

$\hspace{\fill}{\mbox{ et }}
\sum_{i=1}^N\|D_i F\|_H^2\in L^1(\PP^N)\}$

\noindent alors la structure
$$\left(\prod_{n=1}^N(\Omega,\Fc,\PP)_n,\DD_N,\Gamma_N\right)$$ est ferm\'{e}e.

Lorsque $N$ varie elles forment un syst\`{e}me projectif tel que
$$
\begin{array}{rcl}
\Gamma[X_n]&=&\|K_nD_nX_n\|_H^2+\|D_nX_n\|_H^2\\
\Gamma[X_m,X_n]&=&<K_mD_mX_m,K_nD_nX_n>_H\quad\quad m\neq m
\end{array}
$$ qui constitue une pr\'{e}-structure d'erreur, auto-isomorphe
 par translation des indices et non fermable.

Des calculs explicites peuvent \^{e}tre men\'{e}s si l'op\'{e}rateur $K$ est pr\'{e}cis\'{e}. 
 Par exemple si $K$ est la multiplication par $a(\varphi, \omega)$ pour l'espace de Wiener et par $b(\varphi,\omega)$ pour l'angle $\varphi$,

\noindent si $a(\varphi,\omega)=b(\varphi,\omega)\sqrt{g}(\varphi)=1_{[0,2\pi]}(\varphi)$ on obtient
pour l'erreur asymptotique
$$\lim_{N\uparrow\infty}\Gamma[\frac{1}{N}\sum_{n=1}^N A_n]=\int_{[0,L]}\left(\frac{1}{2\pi}
(\int_{\RR}|\cos x|-|\sin x|)\int_t^L\frac{e^{-\frac{x^2}{2s}}}{\sqrt{2\pi s}}\,dsdx\right)^2(dt+\delta_0(t))$$

\noindent et si $a(\varphi,\omega)=b(\varphi,\omega)\sqrt{g}(\varphi)=\varphi.1_{[0,2\pi]}(\varphi)$ on obtient
$$\lim_{N\uparrow\infty}\Gamma[\frac{1}{N}\sum_{n=1}^N A_n]=\int_{[0,L]}\left(\frac{2}{\pi}(L-t)-
\int_{\RR}|\cos x|\int_t^L\frac{e^{-\frac{x^2}{2s}}}{\sqrt{2\pi s}}\,dsdx\right)^2(dt+\delta_0(t))$$

\noindent{\it Remarque}.  Une grande vari\'{e}t\'{e} de mod\`{e}les peuvent \^{e}tre construits et
 manipul\'{e}s avec des hypoth\`{e}ses diff\'{e}rentes sur la ficelle ou sur les erreurs. Ainsi avec une ficelle de classe $\Cc^2$
$$
\begin{array}{rcl}
X(t)&=&X_0+\int_0^t\cos(\varphi+V_s)\,ds\\
Y(t)&=&X_0+\int_0^t\sin(\varphi+V_s)\,ds\quad\quad0\leq t\leq L\leq 1
\end{array}
$$ o\`{u} $V_s$ est le processus gaussien de classe $\Cc^1$
$$V_s=\int_0^1u\wedge s\,dB_u$$ on peut envisager des erreurs qui perturbent la
 ficelle de fa\c{c}on plus profonde (second quantization) par exemple avec la forme de Dirichlet associ\'{e}e au 
semigroupe 
$$(P_tF)(\omega)=\widetilde{\EE}[F(p_t\,\omega+\sqrt{I-p_{2t}}\;\tilde{\omega})]$$ o\`{u}
$p_t$ est le semigroupe de la chaleur sur $[0,1]$ avec r\'{e}flexion au bord,
 la quantit\'{e} mesur\'{e}e
$$A=\int_0^L|\cos(\varphi+V_s)|\,ds$$ restant alors parmi celles qui ont une erreur bien
 d\'{e}finie. Cette structure sur l'espace de Wiener a un  op\'{e}rateur $\Gamma_1$ tel que
$$\Gamma_1[\int_0^1 f(s)\,dB_s]=\int_0^1 f^{\prime 2}(s)\,ds\quad\quad\forall f\in H^1[0,1]$$
formule qui permet comme pr\'{e}c\'{e}demment de le d\'{e}finir de proche en proche par calcul fonctionnel et continuit\'{e}.

On a ici $\Gamma_1[V_s,V_t]=s\wedge t$ d'o\`{u}
$$\Gamma_1[A]=\int_0^L\int_0^L{\mbox{sign}}[\cos(\varphi+V_s)]\sin(\varphi+V_s)
{\mbox{sign}}[\cos(\varphi+V_t)]\sin(\varphi+V_t)(s\wedge t)\,dsdt.$$

\vspace{.5cm}

{\bf Conjectures}\\

La c\'{e}l\`{e}bre m\'{e}thode d'int\'{e}gration par partie en dimension infinie qui permit \`{a} P. Malliavin d'am\'{e}liorer le
 th\'{e}or\`{e}me de H\"{o}rmander sur les op\'{e}rateurs diff\'{e}rentiels du second ordre hypoelliptiques consiste sch\'{e}matiquement
\`{a} d\'{e}finir sur l'espace de Wiener un op\'{e}rateur gradient $\nabla$ dont le transpos\'{e} \'{e}tende l'int\'{e}grale stochastique 
de It\^{o} et \`{a} \'{e}crire un analogue de la formule classique
$$\int g\,{\mbox{div}}\vec{f}=\int \nabla g.\vec{f}$$
$g$ n'\'{e}tant plus ${\cal C}^\infty$ \`{a} support compact mais dans une classe de fonctions tests d\'{e}duite des int\'{e}grales de Wiener.
Cette m\'{e}thode a donn\'{e} plusieurs r\'{e}sultats nouveaux d'existence et de r\'{e}gularit\'{e} ${\cal C}^\infty$
 de densit\'{e} cf [8] et se relie \`{a} la th\'{e}orie des distributions sur l'espace de Wiener cf [13].

Consid\'{e}rer ces op\'{e}rateurs est \'{e}quivalent \`{a} consid\'{e}rer la structure d'erreur d'Ornstein-Uhlenbeck 
que nous avons utilis\'{e}e dans l'exemple C. Elle a la propri\'{e}t\'{e} int\'{e}ressante que toute variable al\'{e}atoire $X$ \`{a} valeurs
dans $\RR^d$ dont la matrice d'erreur $\Gamma[X,X^t]$ est presque s\^{u}rement inversible a une loi absolument continue, ceci peut \^{e}tre d\'{e}montr\'{e} gr\^{a}ce
\`{a} la formule de co-aire de H. Federer. Ceci donne des r\'{e}sultats de r\'{e}gularit\'{e} de solutions d'\'{e}quations diff\'{e}rentielles
stochastiques dans le cas de coefficients lipschitziens. Cette propri\'{e}t\'{e} est vraie pour plusieurs autres structures obtenues par 
image ou produit et est toujours vraie pour des fonctions \`{a} valeurs r\'{e}elles cf [14]. Nous avons pos\'{e} en 1986 avec Francis Hirsch
 la conjecture qu'elle \'{e}tait toujours v\'{e}rifi\'{e}e.\\

\noindent{\bf Conjecture} {\it (conjecture de r\'{e}gularit\'{e} des lois)

Pour toute structure d'erreur $(\Omega, {\cal A}, \PP, \DD, \Gamma)$ et tout entier $d\geq 1$, si $X\in\DD^d$
$$X_{\ast}(1_{\{\det\Gamma[X,X^t]>0\}}.\PP)<<\lambda_d$$
o\`{u} $\lambda_d$ est la mesure de Lebesgue sur $\RR^d$ et $X_{\ast}\mu$ d\'{e}signe l'image de la mesure $\mu$ par $X$.}\\

D\'{e}j\`{a} sur $\RR^d$ si les composantes de l'application identit\'{e} sont dans $\DD$ cette conjecture se relie \`{a} la connaissance
des formes quadratiques diff\'{e}rentielles qui sont fermables qui sont malheureusement encore mal caract\'{e}ris\'{e}es :\\

\noindent{\bf Conjecture} {\it (conjecture des formes ferm\'{e}es)

\it Soit $\mu$ une probabilit\'{e} sur $\RR^n$. Soit $A(x)$ une famille mesurable de matrices  $n\times n$ sym\'{e}triques
telle que la forme d\'{e}finie pour  $f\in{\cal D}(\RR^n)$ par $$\int \nabla f\,A\,(\nabla f)^t\,d\mu$$
soit fermable dans $L^2(\mu)$, alors pour $1\leq k\leq n$
$$1_{\{{\small\rm\mbox{rang de }}A = k\}}.\mu<<{\cal H}^k$$
o\`{u} ${\cal H}^k$ est la mesure de  Hausdorff $k$-dimensionnelle sur $\RR^n$.}\\

Cette propri\'{e}t\'{e} a \'{e}t\'{e} prouv\'{e}e dans le cas  $n=1$ par M. Hamza en 1975 cf [6]. Certains r\'{e}sultats partiellement
non publi\'{e}s de  D. Preiss et 
G. Mokobodzki donnent le cas  $n=k=2$. On attend confirmation.\\

Qu'en est-il du calcul sur les erreurs moyennes dont nous avons dit qu'il fait intervenir les d\'{e}riv\'{e}es secondes ? 
Il n'est pas n\'{e}cessaire d'introduire de nouveaux objets pour le traiter. Comme une forme de Dirichlet sur un espace $L^2(m)$ d\'{e}termine
son semi-groupe d'op\'{e}rateurs et son g\'{e}n\'{e}rateur, l'op\'{e}rateur $L$ v\'{e}rifiant (4) est d\'{e}termin\'{e}, y compris ses termes d'ordre 1, par la structure
$(\Omega, {\cal A}, \PP,\DD,\Gamma)$, simplement il op\`{e}re sur des fonctions plus 
r\'{e}guli\`{e}res et le calcul sur les biais utilise d\'{e}j\`{a} le calcul sur les variances.\\

Mentionnons enfin que d\`{e}s que l'espace $\Omega$ est muni d'une structure topo\-logique qui se relie \`{a} la tribu ${\cal A}$,
une structure d'erreur d\'{e}finit naturellement une notion de capacit\'{e} cf [15] qui ouvre le ``programme quasi-s\^{u}r" de raffinement
 des r\'{e}sultats probabilistes presque s\^{u}rs comme cela a \'{e}t\'{e} fait sur l'espace de Wiener avec
 la structure d'Ornstein-Uhlenbeck.

\vspace{.3cm}

\begin{center}
\noindent {\bf Bibliographie}\\
\end{center}

[1] H. A. Klein, {\it The science of measurement, a historical survey}, Dover, 1974.

[2] K. It\^{o}, ``Stochastic differential equations on a manifold" {\it Nagoya Math. J.} 1, 35-47, 1950.~--- P. A. Meyer, J. Azema, M. 
Yor, (eds) {\it S\'{e}minaire de probabilit\'{e}s XVI, suppl\'{e}ment g\'{e}om\'{e}trie diff\'{e}rentielle stochastique} Lect. Notes in
Math. 921, Springer, 1982~--- N. Ikeda, S. Watanabe, {\it Stochastic
differential equations and diffusion processes} North Holland, 1989~--- M. Emery, {\it Stochastic calculus in manifolds}, Springer, 1989~--- D. Stroock, S. Tanigushi, ``Diffusions as integral curves of vector
fields, Stratonovitch diffusion without It\^{o} integration" in {\it Prog. in Prob.} 34, 333-369, Birk\"{a}user, 1994.

[3] A. Beurling, J. Deny, ``Espaces de Dirichlet, I. le cas \'{e}l\'{e}mentaire" {\it Acta Math.} 99 (1958),
 203-224; ``Dirichlet spaces",
{\it Proc. Nat. Acad.
 Sci. U.S.A.} 45 (1959), 206-215~--- M. L. Silverstein, {\it Symmetric Markov processes},
 L. N. in Math. vol 426, Springer, 1974~--- M. Fukushima, {\it Dirichlet forms and Markov processes},
North-Holland-Kodansha, 1980.

[4] Cl. Dellacherie , P. A. Meyer, {\it Probabilit\'{e}s et potentiel}, Chap. XV \S 2 Hermann 1987~---
N. Bouleau, F. Hirsch, {\it Dirichlet forms and analysis on Wiener space}, Chap.I \S 4, de Gruyter, 1991.

[5] Y. Le Jan, ``Mesures associ\'{e}es \`{a} une forme de Dirichlet" {\it Bull. de la SMF},
 106, (1978), 61-112~--- M. Fukushima {\it op. cit.}~--- J. Bertoin ``Les processus de Dirichlet en tant qu'espace de
Banach" {\it Stochastics} 18, (1986), 155-168~--- T. Lyons, Z. Weian, ``A crossing estimate for the canonical process on a 
Dirichlet space and a tightness result", {\it Ast\'{e}risque}, $n^0$ 157-158, (1988).

[6] M. Fukushima, Y. Oshima, M. Takeda, {\it Dirichlet forms
and Markov processes}, de Gruyter, 1994~---
 Z. Ma, M. R\"{o}ckner, {\it Dirichlet forms}, Springer, 1992.

[7] N. Bouleau, F. Hirsch, {\it op. cit.} Chap.V~--- N. Bouleau, ``Construction
 of Dirichlet structures", {\it Potential theory ICPT 1994}, de Gruyter, 1995.

[8] D. W. Stroock, ``The Malliavin calculus, a functional analytic approach", {\it J. Functional Analysis} 44, p212-257, 1981 --- 
D. Ocone, ``A guide to the stochastic calculus of variations", in {\it Stochastic Analysis and
Related Topics}, H. Korezlioglu, S. Ustunel, eds, L. N. in Math. vol. 1316 Springer 1987 --- J. Potthoff, ``White noise approach to
 Mal\-liavin calculus", {\it J. Functional Analysis}, 71, p207-217, 1987 --- N. Bouleau, F. Hirsch, {\it ibid.}
 Chap.II, III et IV~--- D. Nualart, {\it 
The Malliavin calculus and related topics}, Springer 1995~--- N. Privault, ``In\'{e}galit\'{e}s de Meyer sur l'espace de Poisson" {\it C. R. 
Acad. Sci. Paris}, sI, t.318, p559, 1994 --- S. Ust\"{u}nel, {\it An introduction to analysis on Wiener space},
 L. N. in Math. 1610, Springer 1995~--- P. Malliavin, {\it Stochastic analysis}, Springer, 1997.

[9] E. Hopf, ``On causality, statistics and probability" {\it J. Math. and Physics} MIT, 13, 51-102, 1934~--- E. Engel, ``A road to randomness
in physical systems" {\it Lect. notes in Stat.} 71, Springer 1992.

[10] O. Enshev, D. Stroock, ``Toward a Riemannian geometry on the path space of a Riemannian manifold" {\it J. of Functional Analysis}
 134, 392-416, 1996~--- A. B. Cruzeiro, P. Malliavin, ``Non perturbative construction of invariant measures through confinement by curvature",
 {\it J. Math. Pures et Appl.} 77, 527-537, 1998.

 [11] D. Dacunha-Castelle, M. Duflo, {\it Probabilit\'{e} et statistiques 1.
 Probl\`{e}mes \`{a} temps fixe}, Masson, 1982.

[12] I. A. Ibragimov, R. Z. Has'minskii, {\it Statistical estimation}, Springer 1981~--- V. Genon-Catalot, D. Picard,
 {\it El\'{e}ments de statistique asymptotique} SMAI-Springer 1993.

[13] I. Shigekawa, ``Derivatives of Wiener functionals and absolute continuity of induced measures",
 {\it J. Math. Kyoto Univ.} 20-2, p263-289, 
1980~--- S. Watanabe, {\it Lectures on Stochastic differential equations and Mal\-liavin calculus} Tata Institute, Springer
 1984. 

[14] N. Bouleau, ``D\'{e}composition de l'\'{e}nergie par niveau de potentiel", {\it L. N. in Math.} Vol 1096, Springer 1984.

[15] A. Beurling, J. Deny, {\it op. cit.}~--- G. Choquet, ``Theory of capacities" {\it Ann. Inst. Fourier}, 5, 1953-54~---
J. Deny, ``Th\'{e}orie de la capacit\'{e} dans les espaces fonctionnels" {\it S\'{e}minaire Brelot-Choquet-Deny}, 9\`{e}me ann\'{e}e, 
1, 1964-65~--- P. Malliavin, ``Implicit functions of finite corank on the Wiener space" {\it Tanigushi Int. Symp. Stoch. An., Katata
1982}, Kinoluniya, Tokyo, 1983~--- D. Feyel, A. de la Pradelle,``Capacit\'{e}s gaussiennes"
 {\it Ann. Inst. Fourier} 41 (1991), 49-76.\\

\noindent Nicolas Bouleau

\noindent Direction de la Recherche

\noindent Ecole des Ponts, Paris

\noindent e-mail:~bouleau@enpc.fr

\end{document}